\documentclass[11pt]{article}

\usepackage{amsmath,amsthm,amsfonts,amssymb,latexsym}
\usepackage{amsmath}

\usepackage{a4wide}

\newtheorem{thm}{Theorem}[section]
\newtheorem{cor}[thm]{Corollary}
\newtheorem{lem}[thm]{Lemma}
\newtheorem{prop}[thm]{Proposition}

\newtheorem{rmk}[thm]{Remark}

\begin{document}


\title{Cohomology of $\mathfrak {osp}(2|2)$ acting
on the spaces of linear differential operators on the superspace
$\mathbb{R}^{1|2}$ }

\author{ Nizar Ben Fraj\and Maha Boujelben\thanks{Institut Sup\'{e}rieur
de Sciences Appliqu\'{e}es et Technologie, Sousse, and
D\'epartement de Math\'ematiques, Facult\'e des Sciences de Sfax,
BP 802, 3038 Sfax, Tunisie. E.mails:
benfraj\_nizar@yahoo.fr,\,Maha.Boujelben@fss.rnu.tn}}

\maketitle

\begin{abstract}
We compute the first differential cohomology of the
orthosymplectic Lie superalgebra $\mathfrak{osp}(2|2)$ with
coefficients in the superspace of linear differential operators
acting on the space of weighted densities on the
(1,\,2)-dimensional real superspace. We also compute the same, but
$\mathfrak{osp}(1|2)$-relative, cohomology. We explicitly give
1-cocycles spanning these cohomology. This work is a simplest
generalization of a result by Basdouri and Ben Ammar [Cohomology
of $\frak {osp}(1|2)$ with coefficients in
$\frak{D}_{\lambda,\mu}$. Lett.
 Math. Phys.{\bf 81}, 239--251 (2007)].
\end{abstract}

\maketitle {\bf Mathematics Subject Classification} (2000). 53D55



\section{Introduction}


The space of weighted densities with weight $\lambda$ (or
$\lambda$-densities) on $\mathbb{R}$, denoted by:
\begin{equation*}
{\mathcal F}_\lambda=\left\{ f(dx)^{\lambda}\mid f\in
C^\infty(\mathbb{R})\right\},\quad \lambda\in\mathbb{R},
\end{equation*}
is the space of sections of the line bundle
$(T^*\mathbb{R})^{\otimes^\lambda}$ for positive integer
$\lambda$. Let ${\rm Vect}(\mathbb{R})$ be the Lie algebra of all
vector fields $X_F=F\frac{d}{dx}$ on $\mathbb{R}$, where $F\in
C^\infty(\mathbb{R})$. The {\it Lie derivative} $L_{D}$ along the
vector field $D$ makes ${\mathcal F}_\lambda$ a ${\rm
Vect}(\mathbb{R})$-module for any $\lambda\in\mathbb{R}$:
\begin{equation}\label{Lie1}
L_{X_F}(f(dx)^{\lambda})=L_{X_F}^\lambda(f)(dx)^{\lambda}\quad
\text{with}\quad L_{X_F}^\lambda(f)=Ff'+\lambda fF',
\end{equation}
where $f'$, $F'$ are $\frac{df}{dx}$, $\frac{dF}{dx}$. On the
space $\mathrm{D}_{\lambda,\mu}$ of differential operators
$\mathcal{F}_\lambda\to \mathcal{F}_\mu$ a ${\rm
Vect}(\mathbb{R})$-module structure is given by the formula:
\begin{equation}\label{Lieder2}{X_F}\cdot A=L_{X_F}^\mu\circ A-A\circ
L_{X_F}^\lambda,
\end{equation}
for any differential operator
$A:f(dx)^\lambda\mapsto(Af)(dx)^\mu$.

Lecomte, in \cite{lec}, found the cohomology
$\mathrm{H}^1_\mathrm{diff}\left(\mathfrak{sl}(2),
\mathrm{D}_{\lambda,\mu}\right)$ and
$\mathrm{H}^2_\mathrm{diff}\left(\mathfrak{sl}(2),
\mathrm{D}_{\lambda,\mu}\right)$, where $\mathfrak{sl}(2)$ is
realized as the Lie subalgebra of $\mathrm{Vect}(\mathbb{R})$
spanned by $\left\{X_1,\,X_{x},\,X_{x^2}\right\}$ and where
$\mathrm{H}^*_\mathrm{diff}$ denotes the differential cohomology;
that is, only cochains given by differential operators are
considered. These spaces appear naturally in the problem of
describing the deformations of the
$\mathfrak{sl}(2)$-module 
${\cal S}_{\mu-\lambda} = \bigoplus_{k=0}^\infty
\mathcal{F}_{\mu-\lambda-k}$, the space of symbols of differential
operators of $\mathrm{D}_{\lambda,\mu}$. More precisely, the
elements of $\mathrm{H}^1\left(\mathfrak{sl}(2),V\right)$ classify
the infinitesimal deformations of a $\mathfrak{sl}(2)$-module $V$
and the obstructions to integrability of a given infinitesimal
deformation of $V$ are elements of
$\mathrm{H}^2\left(\mathfrak{sl}(2),V\right)$ (for examples, see
\cite{aalo, abbo, bb1, nr}).

Now, we can study the corresponding super structures. More
precisely, we consider the superspace $\mathbb{R}^{1|n}$ equipped
with a contact  1-form $\alpha_n$, and introduce the superspace
$\mathfrak{F}^n_\lambda$ of $\lambda$-densities on the superspace
$\mathbb{R}^{1|n}$. The spaces $\mathfrak{F}^n_\lambda$ are
modules over $\mathcal{K}(n)$, the Lie superalgebra of contact
vector fields on $\mathbb{R}^{1|n}$; the space
$\mathfrak{D}^n_{\lambda,\mu}$ of differential operators
$\mathfrak{F}^n_\lambda\to \mathfrak{F}^n_\mu$ is, naturally, a
$\mathcal{K}(n)$-module. The spaces
$\mathrm{H}^i_\mathrm{diff}\left(\mathfrak
{osp}(1|2),\mathfrak{D}^n_{\lambda,\mu}\right) $ for $i=1$ and 2
need to be computed in order to describe deformations of the
$\mathfrak{osp}(1|2)$-module $ {\frak
S}^n_{\mu-\lambda}=\bigoplus_{k\geq0}{\frak
F}^n_{\mu-\lambda-\frac{k}{2}}$, a super analogue of ${\cal
S}_{\mu-\lambda}$, see \cite{gmo}.

In \cite {bb}, Basdouri and Ben Ammar studied this question for
$n=1$. In this case, $\mathfrak{sl}(2)$ is replaced by the Lie
superalgebra $\mathfrak{osp}(1|2)$ naturally realized as a
subalgebra of $\mathcal{K}(1)$.

Since there seems to be no conceptual difference in the setting or
results obtained in the study of the cohomology of
$\mathfrak{osp}(n|2)$ acting on the spaces of  linear differential
operators on the superspace $\mathbb{R}^{1|n}$ for $n$ considered
so far (0, 1 and 2 in this paper), the point is not to treat in
further articles the cases $n=3$ and so on.  The point is that the
behavior and certain properties of the Lie superalgebras
$\mathfrak{osp}(n|2)$ and $\mathcal{K}(n)$ are similar for $n< 4$
(\cite{ff}); the cases $n=0$ and $n=1$ are particularly close.
However, in several questions, the case $n=2$ is exceptional due
to an occasional isomorphism $\mathcal{K}(n)\simeq {\rm
Vect}(\mathbb{R}^{1|1})$ (\cite{gls}), and one never knows {\it a
priori} which type of questions will make a given particular $n$
exceptional. We can expect that the properties of
$\mathfrak{osp}(n|2)$ and $\mathcal{K}(n)$ become uniform only for
$n>6$. So somebody has to perform all the calculations in the hope
to find an interesting result (such, for example, as mentioned in
Subsection \ref{ma}).

In this paper we consider the case $n=2$. That is, we consider the
orthosymplectic Lie superalgebra $\mathfrak{osp}(2|2)$ naturally
realized as a subalgebra of $\mathcal{K}(2)$. We compute here
$\mathrm{H}^1_\mathrm{diff}\left(\mathfrak
{osp}(2|2),\mathfrak{D}^2_{\lambda,\mu}\right) $ and
$\mathrm{H}^1_\mathrm{diff}\left(\mathfrak
{osp}(2|2),\mathfrak{osp}(1|2);\mathfrak{D}^2_{\lambda,\mu}\right)
$. Moreover, we give explicit formulae for all the nontrivial
1-cocycles. These spaces arise in the classification of
infinitesimal deformations of the $\mathfrak{osp}(2|2)$-module $
{\frak S}^2_{\mu-\lambda}=\bigoplus_{k\geq0}{\frak
F}^2_{\mu-\lambda-\frac{k}{2}}$. We hope to be able to describe in
the future all the deformations of this module.

\section{Definitions and Notation}

Recall that $C^\infty({\mathbb{R}}^{1|2})$ consists of elements of
the form:
\begin{equation*}
F(x,\theta_1,\theta_2)=f_0(x) + f_1(x)\theta_1 + f_2(x)\theta_2 +
f_{12}(x)\theta_1\theta_2,
\end{equation*}
where $f_0, f_1, f_2, f_{12}\in C^\infty({\mathbb{R}})$, and where
$x$ is the even indeterminate, $\theta_1$ and $\theta_2$ are odd
 indeterminates, i.e., $\theta_i\theta_j=-\theta_j\theta_i$. Let $|F|$ be
the parity of a homogeneous function $F$. Let
\begin{equation*}
\mathrm{Vect}(\mathbb{R}^{1|2})=\left\{F_0\partial_x+F_1\partial_{1}+F_2\partial_{2}\mid
F_i\in C^\infty(\mathbb{R}^{1|2})\right\},
\end{equation*}
where $\partial_{i}=\frac{\partial}{\partial\theta_i}$. Let
$\mathcal{K}(2)$ be the Lie superalgebra of contact vector fields
on $\mathbb{R}^{1|2}$ :
$$
\mathcal{K}(2)=\big\{X\in\mathrm{Vect}(\mathbb{R}^{1|2})~|~\hbox{there
exists}~F\in C^\infty({\mathbb{R}}^{1|2})~ \hbox{such
that}~\mathfrak{L}_X(\alpha_2)=F\alpha_2\big\},
$$
where $\mathfrak{L}_X$ is the Lie derivative along the vector
field $X$ and
\begin{equation*}
\alpha_2=dx+\theta_1 d\theta_1+\theta_2 d\theta_2.
\end{equation*} Any contact vector field on $\mathbb{R}^{1|2}$ can be
expressed as
\begin{equation*}
X_F=F\partial_x-\frac{1}{2}(-1)^{|F|}\sum_{i=1}^2\overline{\eta}_i(F)\overline{\eta}_i,\;\text{
where }\, F\in C^\infty(\mathbb{R}^{1|2})
\end{equation*}
and $\overline{\eta}_i=\partial_{i}-\theta_i\partial_x.$ The
contact bracket is defined by $[X_F,\,X_G]=X_{\{F,\,G\}}$:
\begin{equation}
\{F,G\}=FG'-F'G-\frac{1}{2}(-1)^{|F|}\sum_{i=1}^2\overline{\eta}_i(F)\cdot\overline{\eta}_i(G).
\end{equation}
The orthosymplectic Lie superalgebra $\mathfrak{osp}(2|2)$ can be
realized as a subalgebra of $\mathcal{K}(2)$:
\begin{equation*}\mathfrak{osp}(2|2)=\text{Span}(X_1,\,X_{x},\,X_{x^2},\,X_{x\theta_1},\,
X_{x\theta_2},\, X_{\theta_1},\, X_{\theta_2},\,
X_{\theta_1\theta_2}).\end{equation*} We easily see that
$\mathfrak{osp}(1|2)$ is subalgebra of $\mathfrak{osp}(2|2)$:
$$
\mathfrak{osp}(1|2)=\text{Span}(X_1,\,X_{x},\,X_{x^2},\,X_{x\theta_1},\,
X_{\theta_1})\simeq
\text{Span}(X_1,\,X_{x},\,X_{x^2},\,X_{x\theta_2},\,
X_{\theta_2}).
$$
We define the space of $\lambda$-densities as
\begin{equation}
\mathfrak{F}^2_\lambda=\left\{F(x,\theta_1,\theta_2)\alpha_2^\lambda\mid
F(x,\theta_1,\theta_2) \in C^\infty(\mathbb{R}^{1|2})\right\}.
\end{equation}
As a vector space, $\mathfrak{F}^2_\lambda$ is isomorphic to
$C^\infty(\mathbb{R}^{1|2})$, but the Lie derivative of the
density $G\alpha_2^\lambda$ along the vector field $X_F$ in
$\mathcal{K}(2)$ is now:
\begin{equation}
\label{superaction}
\mathfrak{L}_{X_F}(G\alpha_2^\lambda)=\mathfrak{L}^{\lambda}_{X_F}(G)\alpha_2^\lambda,
\quad\text{with}~~\mathfrak{L}^{\lambda}_{X_F}(G)=\mathfrak{L}_{X_F}(G)+
\lambda F'G.
\end{equation}

Here, we restrict ourselves to the subalgebra
$\mathfrak{osp}(2|2)$, thus we obtain a one-parameter family of
$\mathfrak{osp}(2|2)$-modules on $C^\infty(\mathbb{R}^{1|2})$
still denoted by $\mathfrak{F}^2_\lambda$. As an
$\mathfrak{osp}(1|2)$-module, we have
\begin{equation}\label{isom}
\mathfrak{F}^2_\lambda\simeq\mathfrak{F}^1_\lambda
\oplus\Pi(\mathfrak{F}^1_{\lambda+{1\over2}})
\end{equation}
where $\Pi$ is the change of parity operator.

Since $-\overline{\eta}_i^2=\partial_x,$ and
$\partial_i=\overline{\eta}_i-\theta_i\overline{\eta}_i^2,$ every
differential operator $A\in\frak{D}^2_{\lambda,\mu}$ can be
expressed in the form
\begin{equation}
\label{diff1}
A(F\alpha_2^\lambda)=\sum_{\ell,m}a_{\ell,m}(x,\theta)
\overline{\eta}_1^\ell\overline{\eta}_2^m(F)\alpha_2^\mu,
\end{equation}
where the coefficients $a_{\ell,m}(x,\theta)$ are arbitrary
 functions.
\begin{prop}
\label{decom} As a $\mathfrak{osp}(1|2)$-module, we have
\begin{equation}
\label{amm}
\frak{D}^2_{\lambda,\mu}\simeq\frak{D}^1_{\lambda,\mu}\oplus
\frak{D}^1_{\lambda+\frac{1}{2},\mu+\frac{1}{2}}\oplus
\Pi\left(\frak{D}^1_{\lambda,\mu+\frac{1}{2}}\oplus
\frak{D}^1_{\lambda+\frac{1}{2},\mu}\right).
\end{equation}
\end{prop}
\begin{proofname}. Any element $F\in C^\infty(\mathbb{R}^{1|2})$ can be
uniquely written as follows: $F=F_1+F_2\theta_2$, where
$\partial_2 F_1=\partial_2 F_2=0$. Therefore, for any
$X_H\in\mathfrak{osp}(1|2)$, we easily chek that
\begin{equation*}{\frak
L}^{\lambda}_{X_H}(F)= {\frak L}^{\lambda}_{X_H}(F_1)+{\frak
L}^{\lambda+{1\over2}}_{X_H}(F_2) \theta_2.
\end{equation*}
Thus,  the following map is an $\mathfrak{osp}(1|2)$-isomorphism:
\begin{equation}\label{ph}\begin{array}{llll}
\Phi_{\lambda}:&\mathfrak{F}^2_\lambda
&\rightarrow\mathfrak{F}^1_\lambda \oplus
\Pi(\mathfrak{F}^1_{\lambda+\frac{1}{2}})\\
&F\alpha_2^{\lambda}&\mapsto\left(F_1 \alpha_{1}^{\lambda},~
\Pi(F_2\alpha_{1}^{\lambda+\frac{1}{2}})\right)
\end{array}\end{equation}
 So, we deduce an
$\mathfrak{osp}(1|2)$-isomorphism:
\begin{equation}\label{Phi}
\begin{array}{lll}\Psi_{\lambda,\mu}:&\frak{D}^1_{\lambda,\mu}\oplus
\frak{D}^1_{\lambda+\frac{1}{2},\mu+\frac{1}{2}}\oplus
\Pi\left(\frak{D}^1_{\lambda,\mu+\frac{1}{2}}\oplus
\frak{D}^1_{\lambda+\frac{1}{2},\mu}\right) \rightarrow
\frak{D}^2_{\lambda,\mu}\\[5pt]&A\mapsto\Phi_{\mu}^{-1}\circ
A\circ\Phi_{\lambda}.
\end{array}
\end{equation}Here, we identify
the $\mathfrak{osp}(1|2)$-modules via the following isomorphisms:
\begin{gather*}\begin{array}{llllllll}
\Pi\left(\frak{D}^1_{\lambda,\mu+\frac{1}{2}}\right)&\rightarrow&
\mathrm{Hom_{diff}}\left(\mathfrak{F}^1_\lambda,\Pi(\mathfrak{F}^1_{\mu+\frac{1}{2}})\right)
\quad &\Pi(A)&\mapsto&\Pi\circ A,\\[10pt]
\Pi\left(\frak{D}^1_{\lambda+\frac{1}{2},\mu}\right)&\rightarrow&
\mathrm{Hom_{diff}}\left(\mathfrak{F}^1_{\lambda+\frac{1}{2}},\Pi(\mathfrak{F}^1_{\mu})\right)
\quad &\Pi(A)&\mapsto& A\circ\Pi,\\[10pt]
\frak{D}^1_{\lambda+\frac{1}{2},\mu+\frac{1}{2}}&\rightarrow&
\mathrm{Hom_{diff}}\left(\Pi(\mathfrak{F}^1_{\lambda+\frac{1}{2}}),
\Pi(\mathfrak{F}^1_{\mu+\frac{1}{2}})\right)
\quad &\Pi(A)&\mapsto&\Pi\circ A\circ\Pi.\\[10pt]
\end{array}
\end{gather*}\hfill$\Box$
\end{proofname}
\section{The space $\mathrm{H}^1(\mathfrak {osp}(2|2),\mathfrak{D}^2_{\lambda,\mu})$}
\subsection{Lie superalgebra cohomology, see \cite{c} }
Let $\frak{g}=\frak{g}_{\bar 0}\oplus \frak{g}_{\bar 1}$ be a Lie
superalgebra acting on a superspace $V=V_{\bar 0}\oplus V_{\bar
1}$ and let $\mathfrak{k}$ be a subalgebra of $\mathfrak{g}$. (If
$\frak{k}$ is omitted it assumed to be $\{0\}$.) The space of
$\frak k$-relative $n$-cochains of $\frak{g}$ with values in $V$
is the $\frak{g}$-module
\begin{equation*}
C^n(\frak{g},\frak{k}; V ) := \mathrm{Hom}_{\frak
k}(\Lambda^n(\frak{g}/\frak{k});V).
\end{equation*}
The {\it coboundary operator} $ \delta_n: C^n(\frak{g},\frak{k}; V
)\rightarrow C^{n+1}(\frak{g},\frak{k}; V )$ is a $\frak{g}$-map
satisfying $\delta_n\circ\delta_{n-1}=0$. The kernel of
$\delta_n$, denoted $Z^n(\mathfrak{g},\frak{k};V)$, is the space
of $\frak k$-relative $n$-{\it cocycles}, among them, the elements
in the range of $\delta_{n-1}$ are called $\frak k$-relative
$n$-{\it coboundaries}. We denote $B^n(\mathfrak{g},\frak{k};V)$
the space of $n$-coboundaries.

By definition, the $n^{th}$ $\frak k$-relative  cohomolgy space is
the quotient space
\begin{equation*}
\mathrm{H}^n
(\mathfrak{g},\frak{k};V)=Z^n(\mathfrak{g},\frak{k};V)/B^n(\mathfrak{g},\frak{k};V).
\end{equation*}
We will only need the formula of $\delta_n$ (which will be simply
denoted $\delta$) in degrees 0 and 1: for $v \in
C^0(\frak{g},\,\frak{k}; V) =V^{\frak k}$,~ $\delta v(g) : =
(-1)^{|g||v|}g\cdot v$, where
\begin{equation*}
V^{\frak k}=\{v\in V~\mid~h\cdot v=0\quad\text{ for all }
h\in\frak k\},
\end{equation*}
and  for  $ \Upsilon\in C^1(\frak{g}, \frak{k};V )$,
\begin{equation*}\delta(\Upsilon)(g,h):=
(-1)^{|g||\Upsilon|}g
\cdot\Upsilon(h)-(-1)^{|h|(|g|+|\Upsilon|)}h\cdot
\Upsilon(g)-\Upsilon([g,h])\quad\text{for any}\quad g,h\in
\frak{g}.\end{equation*}

\subsection{The main theorem}
The main result in this paper is the following:

\begin{thm}
\label{main} The space
$\mathrm{H^1_{diff}}(\mathfrak{osp}(2|2),\frak{D}^2_{\lambda,\mu})$
is purely even. It has the following structure:
\begin{equation}
\mathrm{H^1_{diff}}({\mathcal
K}(2),\frak{D}^2_{\lambda,\mu})\simeq\left\{
\begin{array}{llllll}
\mathbb{R}^2&\text{if}\quad \mu-\lambda=0,\\[2pt]
\mathbb{R}^3&\text{if}\quad
(\lambda,\mu)=(-\frac{k}{2},\frac{k}{2}) \hbox{ and }
k\in\mathbb{N}\backslash\{ 0\},\\[2pt]
0&\text{otherwise}.
\end{array}
\right.
\end{equation}

The following 1-cocycles span the corresponding cohomology spaces:

 \begin{equation}\label{maincocyc}
  \begin{array}{lllllllllll}{\Upsilon}_{\lambda,\lambda}(X_G)&=&G'\\[5pt]
\widetilde{\Upsilon}_{\lambda,\lambda}(X_G)&=&\left\{\begin{array}{ll}
\overline{\eta}_1\overline{\eta}_2(G)\hfill\text{
 if }\lambda=0\\[2pt]2\,\lambda\,
\overline{\eta}_1\left(\partial_{2}(G)\right)-(-1)^{|G|}
\left(\partial_{2}(G)\overline{\eta}_1+\theta_2\overline{\eta}_2
\overline{\eta}_{1}(G)\overline{\eta}_2 \right)\hfill\text{ if
}\lambda\neq0\end{array}\right. \\[10pt]
\Upsilon_{-\frac{k}{2},\frac{k}{2}}(X_G)&=&G'\overline{\eta}_1\overline{\eta}_2^{2k-1}\\[5pt]
\widetilde{\Upsilon}_{-\frac{k}{2},\frac{k}{2}}(X_G)&=&
k\overline{\eta}_1(\partial_{2}(G))\overline{\eta}_1\overline{\eta}_2^{2k-1}
-(-1)^{|G|}\left(\partial_{2}(G) \overline{\eta}_{2}^{2k+1}-
\overline{\eta}_1(\theta_2\partial_{2}(G))\overline{\eta}_1^{2k+1}\right)
\\[5pt]
\overline{\Upsilon}_{-\frac{k}{2},\frac{k}{2}}(X_G)&=&(k-1)G''
\overline{\eta}_1\overline{\eta}_2^{2k-3}+ (-1)^{|G|}\left(
\overline{\eta}_2(G')\overline{\eta}_1^{2k-1}-\overline{\eta}_1(G')
\overline{\eta}_2^{2k-1}\right)\end{array}
\end{equation}
\end{thm}
The proof of Theorem \ref{main} will be the subject of Section
\ref{proof}. In fact, we need first the description of
$\mathrm{H}^1_\mathrm{diff}(\mathfrak{osp}(1|2),
\mathfrak{D}^2_{\lambda,\mu})$ and the
$\mathfrak{osp}(1|2)$-relative cohomology
$\mathrm{H}^1_\mathrm{diff}(\mathfrak{osp}(2|2),\mathfrak{osp}(1|2);
\mathfrak{D}^2_{\lambda,\mu})$. To describe the latter one, we
need also the description of some bilinear
$\mathfrak{osp}(1|2)$-invariant mappings.
\section{Invariant Operators and Cohomology of $\frak {osp}(1|2)$}
\subsection{Invariant bilinear differential operators}

Observe that, as a $\frak {osp}(1|2)$-module, we have
$$
\frak{osp}(2|2)\simeq \frak {osp}(1|2) \oplus \Pi(\frak {h}),
$$
where $\frak {h}$ is the subspace of $\mathfrak{F}^1_{-{1\over2}}$
spanned by
$\{\theta_1\alpha_{1}^{-{1\over2}},\,x\alpha_{1}^{-{1\over2}},\,\alpha_{1}^{-{1\over2}}\}$.
In fact, it is easy to see that, for the adjoint action, the Lie
superalgebra $\mathcal {K}(2)$ is isomorphic to
$\mathfrak{F}^2_{-1}$ which is isomorphic, as $\frak
{osp}(1|2)$-module, to $\mathfrak{F}^1_{-1} \oplus
\Pi(\mathfrak{F}^1_{-\frac{1}{2}})$. So, the space
$\frak{osp}(2|2)$ is isomorphic, as a $\frak{osp}(1|2)$-module, to
$\Phi_\lambda(\frak{osp}(2|2))$, where $\Phi_\lambda$ is given by
(\ref{ph}). More precisely, any element $X_F$ is decomposed into
$X_F=X_{F_1}+X_{ F_2\theta_2}$ where  $\partial_2 F_1=\partial_2
F_2=0$, and then $X_{F_1}\in\frak {osp}(1|2)$ and $X_{
F_2\theta_2}$ is identified to $\Pi({
F_2}\alpha_{1}^{-{1\over2}})\in\Pi(\frak {h})$ and it will be
denoted $X_{\bar{F}_2}$.

To compute the $\mathfrak{osp}(1|2)$-relative cohomology of
$\frak{osp}(2|2)$, we need the description of
$\mathfrak{osp}(1|2)$-invariant mappings form
$\mathfrak{h}\otimes\mathfrak{F}^1_\lambda$ to
$\mathfrak{F}^1_{\mu}$. To do that, we first, describe the
$\mathfrak{sl}(2)$-invariant mappings form
$\mathfrak{h}\otimes{\mathcal F}_\lambda$ to ${\mathcal F}_\mu$.
Obviously, as a $\mathfrak{sl}(2)$-module, we have
$\mathfrak{h}\simeq\mathfrak{h}_0\oplus\Pi(\mathfrak{h}_1),$ where
$\mathfrak{h}_0$ is the subspace of $\mathcal{F}_{-{1\over2}}$
spanned by $\{x (dx)^{-{1\over2}},\,(dx)^{-{1\over2}}\}$ and
$\mathfrak{h}_1$ is the subspace of $\mathcal{F}_{0}$ spanned by
1.
\begin{lem}\label{inva} (see \cite{bb})
Let $ A:\mathfrak{h}_0\otimes{\mathcal
F}_\lambda\rightarrow{\mathcal F}_\mu, \,(h
dx^{-{1\over2}},fdx^\lambda)\mapsto A(h,f)dx^{\mu} $ be an
$\mathfrak{sl}(2)$-invariant nontrivial bilinear differential
operator. Then $\mu=\lambda+k-\frac{1}{2}$ where $k$ is a
non-negative integer satisfying
$$k(k-1)(2\lambda+k-1)(2\lambda+k-2)=0,$$ and, up to a scalar
factor, the map $A$ is given by:
\begin{equation*}
A(h,f)=hf^{(k)}+k(2\lambda+k-1)h'f^{(k-1)}.
\end{equation*}
\end{lem}
By a straightforward computation, we can also check the following
lemma.
\begin{lem}\label{invb}
Let $ B:\mathfrak{h}_1\otimes{\mathcal
F}_\lambda\rightarrow{\mathcal F}_\mu, \,(h ,fdx^\lambda)\mapsto
B(h,f)dx^{\mu} $ be a nontrivial $\mathfrak{sl}(2)$-invariant
bilinear differential operator, then
\begin{equation*}\begin{array}{l}
\mu=\lambda\quad\text{or}\quad(\lambda,\,\mu)=({1-k\over2},\,{1+k\over2})
\quad\text{and}\quad B(h,f)= ahf^{(\mu-\lambda)},\end{array}
\end{equation*}where
$k\in\mathbb{N}$ and $a\in\mathbb{R}$.
\end{lem}
\begin{prop}\label{inv}
Let $
\mathcal{A}:\mathfrak{h}\times\mathfrak{F}^1_\lambda\rightarrow\mathfrak{F}^1_\mu,
\,(H\alpha_{1}^{-{1\over2}},F\alpha_{1}^\lambda)\mapsto
\mathcal{A}(H,F)\alpha_{1}^{\mu}$ be a non-zero
$\mathfrak{osp}(1|2)$-invariant bilinear differential operator.
Then one of the following holds:
\begin{itemize}
  \item [i)] If $\mu=\lambda+k-\frac{1}{2}$ where $k$ is a non-negative integer satisfying
 $k(k-1)(2\lambda+k-1)=0,$ then, up to a scalar factor, the map $\mathcal{A}$ is given by:
\begin{equation}
\label{n1}
\mathcal{A}(H,F)=HF^{(k)}+k(2\lambda+k-1)H'F^{(k-1)}-(-1)^{|H|}k
\overline{\eta}_1(H)\overline{\eta}_1(F^{(k-1)}).
\end{equation}
\item [ii)] If $\mu=\lambda+k$, where $k$ is a non-negative integer satisfying
$k(2\lambda+k)(2\lambda+k-1)=0,$ then, up to a scalar factor, the
map $\mathcal{A}$ is given by:
\begin{equation}
\label{n2}
\mathcal{A}(H,F)=(-1)^{|H|}H\overline{\eta}_1(F^{(k)})+(2\lambda+k)\left(\overline{\eta}_1(H)F^{(k)}+
k H'\overline{\eta}_1(F^{(k-1)})\right).
\end{equation}
\end{itemize}
\end{prop}
\begin{rmk}
For $k=0,1$, the operators (\ref{n1}) and (\ref{n2}) are not only
$\mathfrak{osp}(1|2)$-invariant, but also
$\mathcal{K}(1)$-invariant.
\end{rmk}
\begin{proofname}.
Let $A=A_{\bar0}+A_{\bar1}$ be the decomposition of $A$ into even
and odd parts. As $\mathfrak{sl}(2)$-module, we have
\begin{equation}\label{***}
\mathfrak{h}\times\mathfrak{F}^1_\lambda\simeq\mathfrak{h}_0\otimes{\mathcal
F}_\lambda\oplus\mathfrak{h}_0\otimes\Pi({\mathcal
F}_{\lambda+\frac{1}{2}})\oplus\Pi(\mathfrak{h}_1)\otimes{\mathcal
F}_\lambda\oplus\Pi(\mathfrak{h}_1)\otimes\Pi({\mathcal
F}_{\lambda+\frac{1}{2}}).\end{equation}
 So, the map $A_{\bar0}$ is decomposed into four maps:
\begin{equation}\label{even}\begin{array}{llllllll}
& \mathfrak{h}_0\otimes{\mathcal
F}_\lambda&\longrightarrow&{\mathcal F}_\mu, \quad
&\mathfrak{h}_0\otimes\Pi({\mathcal F}_{\lambda+{1\over2}})&\longrightarrow&\Pi({\mathcal F}_{\mu+{1\over2}}),\\
&\Pi(\mathfrak{h}_1)\otimes{\mathcal
F}_\lambda&\longrightarrow&\Pi({\mathcal F}_{\mu+{1\over2}}),
\quad &\Pi(\mathfrak{h}_1)\otimes\Pi({\mathcal
F}_{\lambda+{1\over2}})&\longrightarrow&{\mathcal F}_\mu
\end{array}\end{equation}
and $A_{\bar1}$ is also decomposed into four maps:
\begin{equation}\label{odd}\begin{array}{llllllll}
& \mathfrak{h}_0\otimes{\mathcal
F}_\lambda&\longrightarrow&\Pi({\mathcal F}_{\mu+{1\over2}}),
\quad
&\mathfrak{h}_0\otimes\Pi({\mathcal F}_{\lambda+{1\over2}})&\longrightarrow&{\mathcal F}_{\mu},\\
&\Pi(\mathfrak{h}_1)\otimes{\mathcal
F}_\lambda&\longrightarrow&{\mathcal F}_{\mu}, \quad
&\Pi(\mathfrak{h}_1)\otimes\Pi({\mathcal
F}_{\lambda+{1\over2}})&\longrightarrow&\Pi({\mathcal
F}_{\mu+{1\over2}}).
\end{array}\end{equation}
Observe that the change of parity $\Pi$ commutes with the
$\mathfrak{sl}(2)$-action, therefore, according to Lemma
\ref{inva} and Lemma \ref{invb}, we can deduce the expressions of
the operators (\ref{even}) and (\ref{odd}). We conclude by using
the invariance property with respect to $X_{\theta_1}$ and
$X_{x\theta_1}$. \hfill$\Box$\end{proofname}

\subsection{The space
$\mathrm{H^1_{diff}}(\mathfrak{osp}(1|2),\frak{D}^2_{\lambda,\mu})$}
Let $\mathfrak{g}=\mathfrak{k}\oplus\mathfrak{p}$ be a Lie
superalgebra, where $\mathfrak{k}$ is a subalgebra and
$\mathfrak{p}$ is a $\mathfrak{k}$-module such that
$[\mathfrak{p},\,\mathfrak{p}]\subset\mathfrak{k}$. Consider a
1-cocycle $\Upsilon\in\mathrm{Z}^1(\mathfrak{g},V)$, where $V$ is
a $\mathfrak{g}$-module. The cocycle relation reads
$$
\Upsilon([g,h])-(-1)^{|g||\Upsilon|}g\cdot\Upsilon(h)+(-1)^{|h|(|g|+|\Upsilon|)}
h\cdot\Upsilon(g)=0,\quad g,h\in\mathfrak{g}.
$$
Denote $\Upsilon_{\mathfrak{k}}=\Upsilon_{|\mathfrak{k}}$ and
$\Upsilon_{\mathfrak{p}}=\Upsilon_{|\mathfrak{p}}$. Obviously,
$\Upsilon_{\mathfrak{k}}$ is a 1-coycle over $\mathfrak{k}$ and if
$\Upsilon_{\mathfrak{k}}=0$ then $\Upsilon_{\mathfrak{p}}$ is
$\mathfrak{k}$-invariant. Thus, the space
$\mathrm{H}^1(\mathfrak{g},V)$ is closely related to the space
$\mathrm{H}^1(\mathfrak{k},V)$. Furthermore,
$\Upsilon_{\mathfrak{k}}$ and $\Upsilon_{\mathfrak{p}}$ subject to
the following equations:
\begin{align}
\label{coc2}
&\Upsilon_{\mathfrak{p}}([h,p])-(-1)^{|h||\Upsilon|}h\cdot\Upsilon_{\mathfrak{p}}(p)+
 (-1)^{|p|(|h|+|\Upsilon|)} p\cdot \Upsilon_{\mathfrak{k}}(h)=0,\,
 h\in{\mathfrak{k}}, \,p\in{\mathfrak{p}},\\[4pt]
¨\label{coc3}
&\Upsilon_{\mathfrak{k}}([p,p'])-(-1)^{|p||\Upsilon|}
p\cdot\Upsilon_{\mathfrak{p}}(p')+(-1)^{|p'|(|p|+|\Upsilon|)}
p'\cdot\Upsilon_{\mathfrak{p}}(p)=0,\,
 p, \,p'\in{\mathfrak{p}}.
\end{align}
In our situation,
$\mathfrak{g}=\mathfrak{osp}(2|2),\,\mathfrak{k}=\mathfrak{osp}(1|2),
\,\mathfrak{p}=\Pi(\mathfrak{h})$ and
$V=\frak{D}^2_{\lambda,\mu}$. Thus, as a first step towards the
proof of Theorem \ref{main}, we shall need to compute
$\mathrm{H}^1(\mathfrak{osp}(1|2),\mathfrak{D}^2_{\lambda,\mu})$.
 According to isomorphism
(\ref{amm}), we can see that the knowledge of
$\mathrm{H}^1_\mathrm{diff}(\mathfrak{osp}(1|2),\mathfrak{D}^1_{\lambda,\mu})$
allows us to compute
$\mathrm{H}^1_\mathrm{diff}(\mathfrak{osp}(1|2),\mathfrak{D}^2_{\lambda,\mu})$:
\begin{equation}\label{h}\begin{array}{ll}\mathrm{H}^1_\mathrm{diff}(\mathfrak{osp}(1|2),
\mathfrak{D}^2_{\lambda,\mu})\simeq&\mathrm{H}^1_\mathrm{diff}(\mathfrak{osp}(1|2),
\mathfrak{D}^{1}_{\lambda,\mu})\oplus\mathrm{H}^1_\mathrm{diff}(\mathfrak{osp}(1|2),
\mathfrak{D}^{1}_{\lambda+\frac{1}{2},\mu+\frac{1}{2}})\oplus\\[8pt]&\mathrm{H}^1_\mathrm{diff}(\mathfrak{osp}(1|2),
\Pi(\mathfrak{D}^{1}_{\lambda,\mu+\frac{1}{2}}))\oplus\mathrm{H}^1_\mathrm{diff}(\mathfrak{osp}(1|2),
\Pi(\mathfrak{D}^{1}_{\lambda+\frac{1}{2},\mu})).\end{array}\end{equation}
Of course, we can deduce the structure of
$\mathrm{H^1_{diff}}(\mathfrak{osp}(1|2),\Pi(\frak{D}^{1}_{\lambda,\mu}))$
from
$\mathrm{H^1_{diff}}(\mathfrak{osp}(1|2),\frak{D}^{1}_{\lambda,\mu})$.
Indeed, to any $\Upsilon\in Z^1_{\rm
diff}(\mathfrak{osp}(1|2),\frak{D}^{1}_{\lambda,\mu})$ corresponds
$\widehat{\Upsilon}\in Z^1_{\rm
diff}(\mathfrak{osp}(1|2),\Pi(\frak{D}^{1}_{\lambda,\mu}))$ where
$\widehat{\Upsilon}(X_G)=\Pi\left(\sigma\circ\Upsilon(X_G)\right)$
with $\sigma(F)=(-1)^{|F|}F$. Obviously, $\Upsilon$ is a
coboundary if and only if $\widehat{\Upsilon}$ is a coboundary.
Thus, we recall the space
$\mathrm{H}^1_\mathrm{diff}(\mathfrak{osp}(1|2),\mathfrak{D}^1_{\lambda,\mu})$
which was
 computed in \cite{bb}:
\begin{equation*}
{\rm H}^1_\mathrm{diff}({\frak
{osp}}(1|2),{\frak{D}}^1_{\lambda,\mu})\simeq\left\{
\begin{array}{ll}
\mathbb{R}&\makebox{ if }~~\lambda=\mu, \\[2pt]
\mathbb{R}^2 & \hbox{ if }~~\lambda=\frac{1-k}{2},~\mu=\frac{k}{2},~k\in\mathbb{N}\setminus\{0\},\\[2pt]
0&\makebox { otherwise. }
\end{array}
\right.
\end{equation*}
A basis for the space
$\mathrm{H}^1_\mathrm{diff}(\mathfrak{osp}(1|2),\,\mathfrak{D}^1_{\lambda,\mu})$
is given by the cohomology classes of the $1$-cocycles
$\Gamma_{\lambda,\mu}$ and $\widetilde{\Gamma}_{\lambda,\mu}$
defined by:
\begin{equation}\label{u0}
\begin{array}{llll}
\Gamma_{\lambda,\lambda}(X_G)&=&G',\\
\Gamma_{{1-k\over2},{k\over2}}(X_G)&=&(-1)^{|G|}\overline{\eta}_1^2(G)\overline{\eta}_1^{2k-1},
\\
\widetilde{\Gamma}_{{1-k\over2},{k\over2}}(X_G)&=&(-1)^{|G|}(k-1)\overline{\eta}_1^4(G)\overline{\eta}_1^{2k-3}
+\overline{\eta}_1^3(G)\overline{\eta}_1^{2k-2}.\end{array}
\end{equation}
\subsection{The space
$\mathrm{H}^1_\mathrm{diff}(\mathfrak{osp}(2|2),\mathfrak{osp}(1|2);
\mathfrak{D}^2_{\lambda,\mu})$}
\label{ma} In this subsection we compute the space
$\mathrm{H}^1_\mathrm{diff}(\mathfrak{osp}(2|2),\mathfrak{osp}(1|2);
\mathfrak{D}^2_{\lambda,\mu})$ and we prove that it is nontrivial
which is not the case for $n=1$:
$\mathrm{H}^1_\mathrm{diff}(\mathfrak{osp}(1|2),\mathfrak{sl}(2);
\mathfrak{D}^1_{\lambda,\mu})=0$, see \cite{bb}. Moreover, the
first author, in \cite{bf}, proved that the space
$\mathrm{H}^1_\mathrm{diff}(\mathcal{K}(2),\mathcal{K}(1);
\mathfrak{D}^2_{\lambda,\mu})$ is nontrivial while the space
$\mathrm{H}^1_\mathrm{diff}(\mathcal{K}(1),\mathrm{Vect}(\mathbb{R});
\mathfrak{D}^1_{\lambda,\mu})$ is trivial, see \cite{bbbbk}.
Hence, the case $n=2$ appears as a special case.

\begin{thm}\label{pro} $\dim\mathrm{H}^1_\mathrm{diff}(\mathfrak{osp}(2|2),\mathfrak{osp}(1|2);
\mathfrak{D}^2_{\lambda,\mu})\leq 1$. It is $1$ only if
$\lambda=\mu\neq0$ or $(\lambda,\mu)=(-{k\over2},{k\over2})$,
where $k\in\mathbb{N}\setminus\{0\}$.  The cohomology classes of
the following 1-cocycles generate the corresponding spaces:
\begin{equation} \label{cocylenontrivial}
\begin{array}{lll}
\widetilde{\Upsilon}_{\lambda,\lambda}(X_G)&=&2\,\lambda\,
\overline{\eta}_1\left(\partial_{2}(G)\right)-(-1)^{|G|}
\left(\partial_{2}(G)\overline{\eta}_1+\theta_2\overline{\eta}_2
\overline{\eta}_{1}(G)\overline{\eta}_2 \right),\\[10pt]
\widetilde{\Upsilon}_{-\frac{k}{2},\frac{k}{2}}(X_G)&=&
k\overline{\eta}_1(\partial_{2}(G))\overline{\eta}_1\overline{\eta}_2^{2k-1}
-(-1)^{|G|}\left(\partial_{2}(G) \overline{\eta}_{2}^{2k+1}-
\overline{\eta}_1(\theta_2\partial_{2}(G))\overline{\eta}_1^{2k+1}\right).
\end{array}
\end{equation}
  \end{thm}

To prove Theorem \ref{pro}, we need the following classical fact:
\begin{lem}
\label{classical} Let $\mathfrak{g}$ be a Lie superalgebra and
$A:U\otimes V\rightarrow W$ a bilinear map, where $U,V$ and $W$
are $\mathfrak{g}$-modules . We consider the following associated
maps$$
\begin{array}{llllll}
&A_1:\Pi(U)\otimes V\rightarrow W,
&A_2:\Pi(U)\otimes\Pi(V)\rightarrow\Pi(W),\\
&A_3:\Pi(U)\otimes V\rightarrow\Pi(W),
&A_4:\Pi(U)\otimes\Pi(V)\rightarrow W
\end{array}
$$ defined by $$
\begin{array}{llllll}
&A_1(\Pi(u)\otimes v)&=&(-1)^{|u|}A(u\otimes v),\quad
A_2(\Pi(u)\otimes \Pi(v))&=&(-1)^{|u|}\Pi\left(A(u\otimes
v)\right),\\
&A_3(\Pi(u)\otimes v)&=&(-1)^{|v|}\Pi\left(A(u\otimes
v)\right),\quad A_4(\Pi(u)\otimes \Pi(v))&=&(-1)^{|v|}A(u\otimes
v).
\end{array}
$$ The maps $A_1,A_2,A_3$ and $A_4$ are $\mathfrak{g}$-invariant
if and only if $A$ is $\mathfrak{g}$-invariant.
\end{lem}

\medskip

\begin{proofname}.(Theorem \ref{pro}): Consider a 1-cocycle
$\Upsilon$ over $\mathfrak{osp}(2|2)$ vanishing on
$\mathfrak{osp}(1|2)$. Thus, the equations (\ref{coc2}) and
(\ref{coc3}) become
\begin{align}
 \label{sltr1} &X_G\cdot\Upsilon(X_{\bar {H}})-(-1)^{|G||\Upsilon|}
 \Upsilon([X_G,X_{\bar {H}}])=0, \\[4pt]
 \label{sltr2}
 &(-1)^{|\bar H_1||\Upsilon|}
{X_{\bar{H}_1}}\cdot\Upsilon(X_{\bar {H}_2})-(-1)^{(|\bar
H_1|+|\Upsilon|)|\bar H_2|}
{X_{\bar{H}_2}}\cdot\Upsilon(X_{\bar{H}_1})=0,
\end{align}
for all $X_{\bar{H}}$, $X_{\bar{H}_1}$,
$X_{\bar{H}_2}\in\Pi(\mathfrak{h})$ and
$X_G\in\mathfrak{osp}(1|2)$. According to the isomorphism
(\ref{ph}), the map ${\Upsilon}$ is decomposed into four
components:
\begin{equation}\label{decomp3}\begin{array}{llllllll}
 &\Pi(\mathfrak{h})\times\mathfrak{F}^1_\lambda&\rightarrow&\mathfrak{F}^1_\mu, \quad
&\Pi(\mathfrak{h})\times\Pi(\mathfrak{F}^1_{\lambda+{1\over2}})
&\rightarrow&\Pi(\mathfrak{F}^1_{\mu+{1\over2}}),\\
&\Pi(\mathfrak{h})
\times\mathfrak{F}^1_\lambda&\rightarrow&\Pi(\mathfrak{F}^1_{\mu+{1\over2}}),
\quad &\Pi(\mathfrak{h})
\times\Pi(\mathfrak{F}^1_{\lambda+{1\over2}})&\rightarrow&\mathfrak{F}^1_\mu.
\end{array}\end{equation}
The equation (\ref{sltr1}) expresses the
$\mathfrak{osp}(1|2)$-invariance of each of these bilinear maps.
 Thus, using Proposition \ref{inv},
Lemma \ref{classical} and equation (\ref{sltr2}), we prove that,
if ${\Upsilon}$ is an odd 1-cocycle then, up to a scalar factor,
${\Upsilon}$ is given by ( with $a, b\in\mathbb{R}$ and
$k\in\mathbb{N}\backslash\{ 0\}$):
\begin{equation*}{\Upsilon}=\left\{\begin{array}{lllll}
\delta\left(a\,\partial^k_{{2}}+ b(\bar{\eta}_{1}+
\theta_{2}\bar{\eta}_{1}\bar{\eta}_{2})\partial^{k-1}_x\right)\quad&\text{
if }
\quad(\lambda,\mu)=(\frac{1-k}{2},\frac{k}{2}),\\[6pt]
\delta\left(a\,\partial^k_{{2}}+
b\theta_2\bar{\eta}_{1}\bar{\eta}_{{2}}\partial^{k-1}_x\right)\quad&\text{
if }
\quad(\lambda,\mu)=(-\frac{k}{2},\frac{k-1}{2}),\\[6pt]
\delta\left(\partial_{2}\right) \quad&\text{ if
}\quad\mu=\lambda+\frac{1}{2}\hbox{ and }\lambda\neq0,-\frac{1}{2},\\[6pt]
\delta\left(\theta_{2}\right) \quad&\text{ if
}\quad\mu=\lambda-\frac{1}{2},\\[6pt]
0 \quad&\text{ otherwise. }
\end{array}\right.
\end{equation*}

Now, if $\Upsilon$ is an even 1-cocycle, by the same arguments as
above, we get:
\begin{equation*}\Upsilon=\left\{\begin{array}{lllll}
a\widetilde{\Upsilon}_{-\frac{k}{2},\frac{k}{2}}+
b\delta\left(\bar{\eta}_1\partial_2\partial^{k-1}_x\right)
 &\text{ if } (\lambda,\mu)=(-\frac{k}{2},\frac{k}{2}),\\[6pt]
 a\,\widetilde{\Upsilon}_{\lambda,\lambda}+b\delta\left(\theta_2\bar{\eta}_2\right)&\text{
if }\quad\lambda=\mu\neq0,\\[6pt]
\delta\left(\theta_{2}(a\,\bar{\eta}_{1}+b\bar{\eta}_2)\right)&\text{
if }\quad\lambda=\mu=0,\\[6pt] 0 \quad&\text{
otherwise }.
\end{array}\right.
\end{equation*}
where $\widetilde{\Upsilon}_{\lambda,\lambda}$ and
$\widetilde{\Upsilon}_{-\frac{k}{2},\frac{k}{2}}$ are those given
in (\ref{cocylenontrivial}). Therefore, in order to complete the
proof of Theorem \ref{pro}, we have to study the cohomology
classes of the 1-cocycles $\widetilde{\Upsilon}_{\lambda,\lambda}$
and $\widetilde{\Upsilon}_{-\frac{k}{2},\frac{k}{2}}$.
\begin{lem}
\label{nontrivial} The maps
$\widetilde{\Upsilon}_{\lambda,\lambda}$ and
$\widetilde{\Upsilon}_{-\frac{k}{2},\frac{k}{2}}$ are nontrivial
${\mathfrak{osp}}(1|2)$-relative 1-cocycles.
\end{lem}
\begin{proofname}. First, we can easily see that, for any even
element $F\in C^\infty({\mathbb{R}}^{1|2})$,
\begin{equation}
\label{even1}
\widetilde{\Upsilon}_{-\frac{k}{2},\frac{k}{2}}(X_{\theta_1\theta_2})
(F\alpha_2^{-\frac{k}{2}})=-k\overline{\eta}_1\overline{\eta}_2^{2k-1}(F)\alpha_2^{\frac{k}{2}}
\end{equation}
Next, assume that there exists an even operator
$A\in\frak{D}^2_{-\frac{k}{2},\frac{k}{2}}$ such that
$\widetilde{\Upsilon}_{-\frac{k}{2},\frac{k}{2}}$ is equal to
$\delta A,$ that is
\begin{equation}
\label{cob}
\widetilde{\Upsilon}_{-\frac{k}{2},\frac{k}{2}}(X_G)=\mathfrak{L}^{\frac{k}{2}}_{X_G}\circ
A- A\circ \mathfrak{L}^{-\frac{k}{2}}_{X_G}.
\end{equation}
The operator $A$ is of the form (\ref{diff1});
the condition (\ref{cob}) implies that its coefficients are
constants (which is equivalent to the fact that  $X_1\cdot
 A=0$). Then, it is now easy to check that the condition (\ref{cob}) has no
solution: using formula (\ref{superaction}), we can see that the
expression (\ref{even1}) never appear in the right hand side of
(\ref{cob}). This is a contradiction with our assumption.
Similarly, we prove that the cocycle
$\widetilde{\Upsilon}_{\lambda,\lambda}$ is nontrivial. Lemma
\ref{nontrivial} is proved. Thus we have completed the proof of
Theorem \ref{pro}.
\end{proofname}
\hfill$\Box$
\end{proofname}

\begin{cor}\label{cor1}  Up to a coboundary, any 1-cocycle $\Upsilon\in
Z^1_\mathrm{diff}(\frak{ {osp}}(2|2),\frak{D}^2_{\lambda,\mu})$ is
invariant with respect to the vector field $X_1={\partial_x}$.
That is, the map $\Upsilon$ can be expressed with constant
coefficients.
\end{cor}

\begin{proofname}. The 1-cocycle condition reads:
\begin{equation}
\begin{array}{lll}\label{part1}
X_1\cdot\Upsilon(X_F)-
(-1)^{|F||\Upsilon|}{X_F}\cdot\Upsilon(X_1)-
\Upsilon([X_1,X_F])=0.
\end{array}
\end{equation}
But, from (\ref{u0}) and Theorem \ref{pro}, it follows that, up to
a coboundary,  we have $\Upsilon(X_1)=0$, and therefore the
equation (\ref{part1}) becomes
\begin{equation*}
\begin{array}{lll}\label{}
X_1\cdot\Upsilon(X_F)- \Upsilon([X_1,X_F])=0
\end{array}
\end{equation*}
which is nothing but the invariance property of $\Upsilon$ with
respect to  $X_1$. \hfill$\Box$\end{proofname}

\section{Proof of Theorem \ref{main}}\label{proof}

Consider a 1-cocycle $\Upsilon\in
Z^1_\mathrm{diff}(\mathfrak{{osp}}(2|2),\mathfrak{D}^2_{\lambda,\mu})$.
If   $\Upsilon_{|\mathfrak{{osp}}(1|2)}$ is trivial then the
1-cocycle $\Upsilon$ is completely described by Theorem \ref{pro}.
Thus, assume that $\Upsilon_{|\mathfrak{{osp}}(1|2)}$ is
nontrivial. Of course, up to coboundary, the general form of
$\Upsilon_{|\mathfrak{{osp}}(1|2)}$ is given by (\ref{u0})
together with the isomorphism (\ref{h}) while
$\Upsilon_{|\Pi(\mathfrak{h})}$ can be essentially described by
equation (\ref{coc2}) and Corollary \ref{cor1}. More precisely,
according to (\ref{u0}) and the isomorphism (\ref{h}), the
1-cocycle $\Upsilon$
 can be nontrivial only, a priori, if $\lambda=\mu$, or
$\lambda=\mu\pm{1\over2}$, or
$(\lambda,\mu)=(-{k\over2},{k\over2}),\,({1-k\over2},{k\over2}),\,(-{k\over2},{k-1\over2})
$ where $k\in\mathbb{N}.$ Thus, we have to distinguish all these
cases. Hereafter all $\varepsilon$'s are constants.

{\bf The case where $\lambda=\mu$}

\noindent Considering (\ref{u0}) and the isomorphism (\ref{h}), we
see that there are two subcases:

i) $\lambda=\mu\neq0$. In this case, the map
$\Upsilon_{|\mathfrak{{osp}}(1|2)}$ is, a priori, given by
$$\Upsilon_{|\mathfrak{{osp}}(1|2)}(X_{G_1})(F\alpha_2^\lambda)=
\left(\varepsilon_1 G'_1F_1+\varepsilon_2
G'_1F_2\theta_2\right)\alpha_2^\lambda,$$ where
$F=F_1+F_2\theta_2$, with  $\partial_2 F_1=\partial_2 F_2=0$. By
direct computation, using equations (\ref{coc2})--(\ref{coc3}) and
Corollary \ref{cor1}, we deduce that
$$\varepsilon_1=\varepsilon_2\quad\hbox{ and }\quad\Upsilon_{|\Pi(\mathfrak{h})}(X_{\bar{G}_2})(F\alpha_2^\lambda)=
\varepsilon_1(-1)^{|F|} G'_2F\theta_2\alpha_2^\lambda.$$
 Hence $\Upsilon$ is a
multiple of $\Upsilon_{\lambda,\lambda},$ see (\ref{maincocyc}).

ii) $\lambda=\mu=0$. Here the map
$\Upsilon_{|\mathfrak{{osp}}(1|2)}$ is, a priori, $\hbox{given
by}$
$$
\Upsilon_{|\mathfrak{{osp}}(1|2)}(X_{G_1})(F)=\left(\varepsilon_1
G'_1F_1+\left(\varepsilon_2
G'_1F_2+(-1)^{|F_1|}\left(\varepsilon_3G_1'
\overline{\eta}_1(F_1)+\varepsilon_4\eta_1(G_1')F_1\right)\right)\theta_2\right).
$$
The same arguments, as above, show that
$\varepsilon_1=\varepsilon_2,~\varepsilon_3=0$ and
$$\Upsilon_{|\Pi(\mathfrak{h})}(X_{\bar{G}_2})(F)=
\left(\varepsilon_1 G'_2\theta_2+
\varepsilon_4(-1)^{|G_2|}\overline{\eta}_1(G_2)\right)F.$$ Hence
$\Upsilon$ is linear combination of $\Upsilon_{0,0}$ and
$\widetilde{\Upsilon}_{0,0},$ see (\ref{maincocyc}).

{\bf The case where }$\mu-\lambda=k$ and $2\lambda=-k\neq0$

\noindent In this case,  the map
$\Upsilon_{|\mathfrak{{osp}}(1|2)}$ is, a priori, $\hbox{given
by}$
$$
\begin{array}{lll}
\Upsilon_{|\mathfrak{{osp}}(1|2)}(X_{G_1})(F\alpha_2^{-{k\over2}})&=&\Big((-1)^{|F_1|}\left(\varepsilon_1G_1'
\overline{\eta}_1(F_1^{(k)})+\varepsilon_2\left(kG_1''\overline{\eta}_1^{2k-1}(F_1)+\eta_1(G_1')
\overline{\eta}_1^{2k}(F_1)\right)\right)\theta_2\\
&&+(-1)^{|F_2|}\Big(\varepsilon_3\left((k-1)G_1''\overline{\eta}_1^{2k-3}(F_2)+\eta_1(G_1')
\overline{\eta}_1^{2k-2}(F_2)\right)\\
&&+\varepsilon_4G_1'
\overline{\eta}_1(F_2^{(k-1)})\Big)\Big)\alpha_2^{{k\over2}}.
\end{array}
$$
Again by the same arguments, we prove that we have
$\varepsilon_4=-\varepsilon_1,~\varepsilon_3=\varepsilon_2$ and
$$
\begin{array}{lll}
\Upsilon_{|\Pi(\mathfrak{h})}(X_{\bar{G}_2})(F\alpha_2^{-{k\over2}})&=&\left(\varepsilon_1
G'_2\overline{\eta}_1(F_2^{(k-1)})\theta_2-
\varepsilon_2G_2'\overline{\eta}_1^{2k-1}(F)\right)
\alpha_2^{{k\over2}}.
\end{array}
$$
Hence $\Upsilon$ is linear combination of
$\Upsilon_{-\frac{k}{2},\frac{k}{2}}$ and
$\overline{\Upsilon}_{-\frac{k}{2},\frac{k}{2}},$ see
(\ref{maincocyc}).

For the cases where $\lambda=\mu\pm{1\over2}$, or
$(\lambda,\mu)=({1-k\over2},{k\over2}),\,(-{k\over2},{k-1\over2}),$
the same arguments as before, show that $\Upsilon$ is trivial.
This completes the proof.\hfill $\Box$

\medskip

\noindent {\bf Acknowledgements} We are pleased to  thank Valentin
Ovsienko, Claude Roger and Christian Duval for their interest in
this work. Special thanks are due to Mabrouk Ben Ammar for his
constant interest in this work and a number of suggestions that
have greatly improved this paper.



\end{document}